\documentclass[A4paper,12pt] {article}

\usepackage{amsmath}
\usepackage{amssymb}
\usepackage{geometry}
\usepackage[pdftex]{graphicx}

\newtheorem {thm} {Theorem}[section] 
\newtheorem {cly} {Corollary} [section]
\newtheorem {defn} {Definition}[section]
\newtheorem {nt} {Note}[section]
\newtheorem {rk} {Remark}[section]
\newtheorem {eg} {Example}[section]
\begin{document}
\title {FG- coupled fixed point theorems for various contractions in partially ordered metric spaces}
\author {Prajisha Eacha, Shaini Pulickakunnel}

\date{}
\maketitle 
\section*{Abstract}
In this paper we introduce FG- coupled fixed point, which is a generalization of coupled fixed point for nonlinear mappings in partially ordered complete metric spaces. We discuss existence and uniqueness theorems of FG- coupled fixed points for different contractive mappings. Our theorems generalizes the results of Gnana Bhaskar and Lakshmikantham [T. Gnana Bhaskar, V. Lakshmikantham; Fixed point theorems in partially ordered metric spaces and applications; Nonlinear Analysis 65 (2006) 1379 - 1393]. \\

\hspace{-0.6cm}\textbf{Key words}: FG- coupled fixed point; Coupled fixed point; Mixed monotone property; Partially ordered set.\\
\textbf{MSC (2010)}: 47H10, 54F05
\section{Introduction}
Fixed point theory has many applications in nonlinear analysis. In \cite{4,6,5} the authors presented fixed point theorems in partially ordered metric spaces and their applications. As a generalization of fixed points, in \cite{1} Guo and Lakshmikantham introduced the concept of abstract coupled fixed points for some operators, thereafter Gnana Bhaskar and Lakshmikantham in \cite{2} introduced coupled fixed points and mixed monotone property for contractive mappings on partially ordered metric spaces. They proved interesting coupled fixed point results in \cite{2}. An interesting application of their result is that it can be used to find the solution of  periodic boundary value problem, moreover it guarantees the uniqueness of the solution. Followed by this several authors established new coupled fixed point theorems in partially ordered complete metric spaces and in cone metric spaces. In \cite{3} Sabetghadam,  Masiha and Sanatpour proved generalization of results of Gnana Bhaskar and Lakshmikantham  in cone metric spaces.\\\\ In this paper we introduce a new concept which is a generalization of coupled fixed point and prove existence theorems for contractive mappings in partially ordered metric spaces. Some examples are also discussed to illustrate our results. We recall the basic definitions.
\begin{defn}[\cite{2}]
\normalfont Let $(X, \leq)$ be a partially ordered set and $F: X\times X\rightarrow X$. We say that F has the mixed monotone property if $F(x, y)$ is monotone non decreasing in $x$ and is monotone non increasing in $y$, that is for any $x, y\in X$\\ $x_{1}, x_{2}\in X$, $x_{1}\leq x_{2} \Rightarrow F(x_{1}, y)\leq F(x_{2}, y)$ and \\ $y_{1}, y_{2}\in X$, $y_{1}\leq y_{2} \Rightarrow F(x, y_{1})\geq F(x, y_{2})$.
\end{defn}
\begin{defn}[\cite{2}]
\normalfont We call an element $(x, y)\in X\times X$ a coupled fixed point of the mapping F if $F(x, y)=x,\ F(y, x)=y$.
\end{defn}
\section{Main Results}
\begin{defn}
\normalfont Let $(X,  \leq_{P_{1}})$ and $(Y,  \leq_{P_{2}})$ be two partially ordered sets and $F: X\times Y \rightarrow X$ and $G: Y\times X\rightarrow Y$ be two mappings. An element $(x,  y)\in X\times Y$ is said to be an FG- coupled fixed point if $F(x,  y)=x$ and $G(y,  x)=y$.
\end{defn} 
\begin{nt}
\normalfont If $X=Y$ and $F=G$ then FG- coupled fixed point becomes coupled fixed point.
An element $(x,  y)\in X\times Y$ is FG- coupled fixed point $\Leftrightarrow$ $(y,  x)\in Y\times X$ is GF- coupled fixed point.
\end{nt}
\begin{nt}
\normalfont Let $(X,  d_{X},  \leq_{P_{1}})$ and $(Y,  d_{Y},  \leq_{P_{2}})$ be two partially ordered metric spaces, then we define the partial order $\leq$ and metric $d$ on $X\times Y$ as follows:\\
For all $(x,  y), (u,  v)\in X\times Y$, $(x,  y) \leq \ (u,  v)\Rightarrow x\leq_{P_{1}} \ u \ \ and \ \ y\geq_{P_{2}} \ v$ and\\ $d ((x,  y),  (u,  v))= d_{X}(x, u)+ d_{Y}(y,  v)$.
\end{nt}
\begin{defn}\label{defn 1}
\normalfont Let $(X,  \leq_{P_{1}})$ and $(Y,  \leq_{P_{2}})$ be two partially ordered sets and $F: X\times Y \rightarrow X$ and $G: Y\times X\rightarrow Y$. We say that F and G have mixed monotone property if F and G are monotone increasing in first variable and monotone decreasing in second variable, i.e, if for all $(x,  y)\in X\times Y$,\\$x_{1}, \ x_{2}\in X$, \ $x_{1}\leq_{P_{1}} x_{2} \Rightarrow F(x_{1},  y)\leq_{P_{1}} F(x_{2},  y)$ and $G(y,  x_{1})\geq_{P_{2}} G(y,  x_{2})$ and \\$y_{1}, \ y_{2}\in Y$, \ $y_{1}\leq_{P_{2}} y_{2}\Rightarrow F(x,  y_{1}) \geq_{P_{1}} F(x, y_{2})$ and $G(y_{1},  x)\leq_{P_{2}} G(y_{2},  x)$.
\end{defn}
\begin{nt}
\normalfont Let $F:X\times Y \rightarrow X$ and $G: Y\times X\rightarrow Y$ be two mappings, then for $n\geq 1$, $F^{n}(x, y)= F(F^{n-1}(x, y),G^{n-1}(y, x))$ and $G^{n}(y, x)=G(G^{n-1}(y, x),F^{n-1}(x, y))$  where $F^{0}(x, y)= x$ and $G^{0}(y, x)=y$ for all $x\in X$ and $y\in Y$.
\end{nt}

\begin{thm} \label{thm 1}
Let $(X,  d_{X},  \leq_{P_{1}})$ and $(Y,  d_{Y},  \leq_{P_{2}})$ be two partially ordered complete metric spaces and $F: X\times Y \rightarrow X$ and $G: Y\times X\rightarrow Y$ be two continuous functions having the mixed monotone property. Assume that there exist $k, \ l\in [0, 1)$ with 
\begin{equation}\label{eqn 1}
d_{X}(F(x, y), F(u, v))\leq \frac{k}{2} \ [d_{X}(x, u)+ d_{Y}(y, v)], \ \forall \   x\geq_{P_{1}}u,\  \ y\leq_{P_{2}}v 
\end{equation}
\begin{equation}\label{eqn 2}
d_{Y}(G(y, x), G(v, u))\leq \frac{l}{2} \ [d_{Y}(y,  v)+ d_{X}(x, u)], \ \forall  \  x\leq_{P_{1}}u,\  \ y\geq_{P_{2}}v
\end{equation}
If there exist $(x_{0}, y_{0})\in X\times Y$ such that $x_{0}\leq_{P_{1}} F(x_{0}, y_{0})$ and $y_{0}\geq_{P_{2}}G(y_{0}, x_{0})$,
then there exist $(x, y)\in X\times Y$ such that $x= F(x, y)$ and $y= G(y, x)$.
\end{thm}
 \textbf{Proof}: 
By hypothesis there exists $(x_{0}, y_{0})\in X\times Y$ such that \\
$x_{0}\leq_{P_{1}} F(x_{0}, y_{0})= x_{1}$ (say) and $y_{0}\geq_{P_{2}}G(y_{0}, x_{0})= y_{1}$ (say).\\ 
For $n= 1, 2, 3,...$ we define $x_{n+1}= F(x_{n}, y_{n})$ and $ y_{n+1}=G(y_{n}, x_{n})$ then we get\\
$x_{n+1}= F^{n+1}(x_{0}, y_{0})$ and $y_{n+1}= G^{n+1}(y_{0},  x_{0})$.\\  
Then we can easily prove that $\{x_{n}\}$ is an increasing sequence in X and $\{y_{n}\}$ is a decreasing sequence in Y by using the mixed monotone property of F and G.\\ Claim: For $n \in \mathbb{N}$
\begin{equation} \label{eqn 4}
d_{X}(F^{n+1}(x_{0}, y_{0}), F^{n}(x_{0}, y_{0}))\leq \frac{k}{2} \ \Big(\frac{k+l}{2}\Big)^{n-1}[d_{X}(x_{1}, x_{0})+ d_{Y}(y_{1}, y_{0})]
\end{equation}
\begin{equation} \label{eqn 5}
d_{Y}(G^{n+1}(y_{0}, x_{0}), G^{n}(y_{0}, x_{0}))\leq \frac{l}{2} \ \Big(\frac{k+l}{2}\Big)^{n-1}[d_{Y}(y_{1},  y_{0})+ d_{X}(x_{1}, x_{0})]
\end{equation}
We will use the fact that $\{x_{n}\}$ is an increasing sequence in X and $\{y_{n}\}$ is a decreasing sequence in Y, (\ref{eqn 1}), (\ref{eqn 2}) and symmetric property of $d_{Y}$ to prove the claim.\\For $n= 1$,\\
$d_{X}(F^{2}(x_{0}, y_{0}), F(x_{0}, y_{0}))= d_{X}(F (F(x_{0}, y_{0}), G(y_{0}, x_{0})), F(x_{0}, y_{0}))\\ \hspace*{4.4cm}\leq \dfrac{k}{2}\ [d_{X}(F(x_{0}, y_{0}), x_{0})+ d_{Y}(G(y_{0}, x_{0}), y_{0})]\\ \hspace*{4.4cm}= \dfrac{k}{2}\ [d_{X}(x_{1}, x_{0})+ d_{Y}(y_{1}, y_{0})]$\\Similarly        
$d_{Y}(G^{2}(y_{0}, x_{0}), G(y_{0},x_{0})) \leq \dfrac{l}{2}\ [d_{Y}(y_{0}, y_{1})+ d_{X}(x_{0}, x_{1})]$\\ 
Now assume the claim for $n\leq m$ and check for $n= m+ 1$.\\Consider,
$d_{X}(F^{m+2}(x_{0}, y_{0}), F^{m+1}(x_{0}, y_{0}))\\\hspace*{.5cm}= d_{X}(F(F^{m+1}(x_{0}, y_{0}), G^{m+1}(y_{0}, x_{0})), F(F^{m}(x_{0}, y_{0}), G^{m}(y_{0}, x_{0})))\\ \hspace*{.5cm}\leq  \dfrac{k}{2}\ [d_{X}(F^{m+1}(x_{0}, y_{0}), F^{m}(x_{0}, y_{0}))+ d_{Y}(G^{m+1}(y_{0}, x_{0}), G^{m}(y_{0}, x_{0}))]\\\hspace*{.5cm} \leq  \dfrac{k}{2} \ \Big\{ \dfrac{k}{2}\ \Big(\dfrac{k+l}{2}\Big)^{m-1} [d_{X}(x_{1}, x_{0})+ d_{Y}(y_{1}, y_{0})] + \dfrac{l}{2}\ \Big(\dfrac{k+l}{2}\Big)^{m-1} [d_{Y}(y_{1}, y_{0})+ 
 d_{X}(x_{1}, x_{0})] \Big\} \\\hspace*{.5cm}
 = \dfrac{k}{2}\ \Big(\dfrac{k+l}{2}\Big)^{m}[d_{X}(x_{1}, x_{0})+d_{Y}(y_{1}, y_{0})]$\\
Similarly we can show that\\
$d_{Y}(G^{m+2}(y_{0}, x_{0}), G^{m+1}(y_{0}, x_{0}))\leq \ \dfrac{l}{2}\ \Big(\dfrac{k+l}{2}\Big)^{m} [d_{X}(x_{1}, x_{0})+d_{Y}(y_{1}, y_{0})]$\\
Thus the claim is true for all $n\in \mathbb{N}$.
Using the result obtained we prove that $\{x_{n}\}$ is a Cauchy sequence in X and $\{y_{n}\}$ is a Cauchy sequence in Y. \\For $m\geq n$ consider,\\
$d_{X}(F^{m}(x_{0}, y_{0}), F^{n}(x_{0}, y_{0}))\\ \hspace*{.5cm}\leq  d_{X}(F^{m}(x_{0}, y_{0}), F^{m-1}(x_{0}, y_{0}))+ \ d_{X}(F^{m-1}(x_{0}, y_{0}), F^{m-2}(x_{0}, y_{0}))\\\hspace*{.5cm}~~~+ . . .  + d_{X}(F^{n+1}(x_{0}, y_{0}), F^{n}(x_{0}, y_{0})) \\ \hspace*{.5cm}\leq  \dfrac{k}{2}\ \Big( \dfrac{k+l}{2}\Big)^{m-2}\ [d_{X}(x_{1},  x_{0})+ d_{Y}(y_{1},  y_{0}) ]+ \dfrac{k}{2}\ \Big( \dfrac{k+l}{2}\Big)^{m-3}\ [d_{X}(x_{1},  x_{0})+ d_{Y}(y_{1}, y_{0})]\\\hspace*{.5cm}~~~ + . . . + \dfrac{k}{2}\ \Big( \dfrac{k+l}{2} \Big)^{n-1}\ [d_{X}(x_{1}, x_{0})+ d_{Y}(y_{1}, y_{0})] \\ \hspace*{.5cm}= \Big[ \dfrac{k}{2} \ \Big( \dfrac{k+l}{2}\Big)^{m-2}+ \dfrac{k}{2} \ \Big( \dfrac{k+l}{2}\Big)^{m-3}+ . . . + \dfrac{k}{2} \ \Big( \dfrac{k+l}{2}\Big)^{n-1}\Big] [d_{X}(x_{1}, x_{0})+ d_{Y}(y_{1}, y_{0})] \\\hspace*{.5cm}\leq  \dfrac{k}{2}\  \Big(\dfrac{\theta^{n-1}}{1-\theta}\Big) [d_{X}(x_{1}, x_{0})+ d_{Y}(y_{1}, y_{0})] \ ;\ \text{ where} \ \ \theta= \dfrac{k+l}{2} < 1    \\ \hspace*{.5cm} \rightarrow  0 \ \ \text{as} \ \ n\rightarrow \ \ \infty $\\
That is $\{F^{n}(x_{0},y_{0})\}_{n= 0}^{\infty}$ is Cauchy sequence in $(X, d_{X})$.\\Similarly we get $\{G^{n}(y_{0},x_{0})\}_{n= 0}^{\infty}$ is a Cauchy sequence in $(Y, d_{Y})$.\\  
Since $(X, d_{X})$ and $(Y, d_{Y})$ are complete metric spaces, we have $\lim_{n \rightarrow \infty}F^{n}(x_{0}, y_{0})= x$ and $\lim_{n \rightarrow \infty}G^{n}(y_{0}, x_{0})= y$ for some $(x, y)\in X\times Y$. Now we can prove that $(x, y)$ is an FG- coupled fixed point by using the continuity of F and G. For that consider,\\ 
$d_{X}(F(x, y), x)=  \lim_{n\rightarrow \infty} d_{X}(F(F^{n}(x_{0}, y_{0}), G^{n}(y_{0}, x_{0})), F^{n}(x_{0}, y_{0}))\\\hspace*{2.5cm} =  \lim_{n\rightarrow \infty} d_{X}(F^{n+1}(x_{0}, y_{0}), F^{n}(x_{0}, y_{0}))\\\hspace*{2.5cm} =  0$\\
That is $F(x, y)= x$.\\In a similar manner we can prove that $G(y, x)= y$. This completes the proof. $\square$

\begin{eg}
Let $X=(-\infty, 0]$ and $Y=[0,\infty)$ with usual order and usual metric. Define $F:X\times Y \rightarrow X$ and $G:Y \times X \rightarrow Y$ as $F(x,y)=\dfrac{x-y}{3}$ and $G(y,x)=\dfrac{y-x}{5}$, then it is easy to check the conditions $(\ref{eqn 1})$ and $(\ref{eqn 2})$ for F and G with $k=\dfrac{2}{3},\ l=\dfrac{2}{5}$. Here $(0,0)$ is the unique FG- coupled fixed point.
\end{eg}
We obtain the result of Gnana Bhaskar and Lakshmikantham \cite{2} as a corollary of our result.
\begin{cly} {\upshape\lbrack 1, Theorem 2.1 \rbrack} \label{cly 1}
Let $(X, \leq)$ be a partially ordered set and suppose there is a metric d on X such that $(X, d)$ is a complete metric space. Let $F: X\times X \rightarrow X$  be a continuous mapping having the mixed monotone property on X. Assume that there exist $k \in [0, 1)$ with 
\begin{equation*}
d(F(x, y), F(u, v))\leq \frac{k}{2} \ [d(x, u)+ d(y, v)], \ \forall \   x\geq u,\  \ y\leq v 
\end{equation*}
If there exists $x_{0}, y_{0}\in X$ such that $x_{0}\leq F(x_{0}, y_{0})$ and $y_{0}\geq F(y_{0}, x_{0})$,
then there exist $x, y \in X$ such that $x= F(x, y)$ and $y= F(y, x)$.
\end{cly}
\textbf{Proof}: Take $X= Y$, $F=G$ and $k= l$ in Theorem \ref{thm 1}, we get the result. $\square$

\begin{rk}
 \normalfont By adding to the hypothesis of Theorem \ref{thm 1} the condition: for every $(x, y), (x_{1}, y_{1})\in X\times Y$ there exists a $(u, v)\in X\times Y$ that is comparable to both $(x, y)$  and  $(x_{1}, y_{1})$, we can obtain a unique FG- coupled fixed point. 
\end{rk} In the following theorem we prove the uniqueness of FG- coupled fixed point using the above condition.

\begin{thm} 
Let $(X,  d_{X},  \leq_{P_{1}})$ and $(Y,  d_{Y},  \leq_{P_{2}})$ be two partially ordered complete metric spaces and $F: X\times Y \rightarrow X$ and $G: Y\times X\rightarrow Y$ be two continuous functions having the mixed monotone property. Assume that for every $(x, y), (x_{1}, y_{1})\in X\times Y$ there exists a $(u, v)\in X\times Y$ that is comparable to both $(x, y) \ and \ (x_{1}, y_{1})$ and there exist $k, \ l\in [0, 1)$ with 
\begin{equation*}
d_{X}(F(x, y), F(u, v))\leq \frac{k}{2} \ [d_{X}(x, u)+ d_{Y}(y, v)], \ \forall \   x\geq_{P_{1}}u,\  \ y\leq_{P_{2}}v \tag{\ref{eqn 1}}
\end{equation*}
\begin{equation*}
d_{Y}(G(y, x), G(v, u))\leq \frac{l}{2} \ [d_{Y}(y,  v)+ d_{X}(x, u)], \ \forall  \  x\leq_{P_{1}}u,\  \ y\geq_{P_{2}}v \tag{\ref{eqn 2}}
\end{equation*}
If there exist $(x_{0}, y_{0})\in X\times Y$ such that $x_{0}\leq_{P_{1}} F(x_{0}, y_{0})$ and $y_{0}\geq_{P_{2}}G(y_{0}, x_{0})$,
then there exist unique $(x, y)\in X\times Y$ such that $x= F(x, y)$ and $y= G(y, x)$.
\end{thm}
\textbf{Proof}: Following as in Theorem \ref{thm 1} we obtain the existence of FG- coupled fixed point. Now we show the uniqueness part. Suppose that $(x^{*}, y^{*})\in X\times Y$ is another FG- coupled fixed point, then we show that $d((x, y), (x^{*}, y^{*}))= 0$,\\where $x= \lim_{n\rightarrow \infty} F^{n}(x_{0}, y_{0})$ and $y= \lim_{n\rightarrow \infty}G^{n}(y_{0}, x_{0})$.\\ 
Claim: For any two points $(x_{1}, y_{1}),\ (x_{2}, y_{2})\in X\times Y$ which are comparable, 
\begin{equation}\label{eqn 6}
d_{X}(F^{n}(x_{1}, y_{1}), F^{n}(x_{2}, y_{2}))= \Big(\frac{k}{2}\Big)^{n} [d_{X}(x_{1}, x_{2})+ d_{Y}(y_{1}, y_{2})]
\end{equation}
\begin{equation}\label{eqn 7}
d_{Y}(G^{n}(y_{1}, x_{1}), G^{n}(y_{2}, x_{2}))= \Big(\frac{l}{2}\Big)^{n} [d_{Y}(y_{1}, y_{2})+ d_{X}(x_{1}, x_{2})]
\end{equation} Without loss of generality assume that $(x_{2},y_{2})\leq(x_{1},y_{1})$.\\
We will use (\ref{eqn 1}), (\ref{eqn 2}) and symmetric property of $d_{Y}$ to prove the claim. \\For $n= 1$ consider,\\
$d_{X}(F(x_{1}, y_{1}), F(x_{2}, y_{2})) \leq \dfrac{k}{2} \ [d_{X}(x_{1}, x_{2})+d_{Y}(y_{1}, y_{2})]\\ 
d_{Y}(G(y_{1}, x_{1}), G(y_{2}, x_{2}))\leq \dfrac{l}{2}\ [d_{Y}(y_{1}, y_{2})+d_{X}(x_{1}, x_{2})]$\\ 
That is our claim is true for $n= 1$.\\ Assume that it is true for $n\leq m$ and check for $n= m+ 1$. \\Consider,\\
$d_{X}(F^{m+1}(x_{1}, y_{1}), F^{m+1}(x_{2}, y_{2}))=  d_{X}(F(F^{m}(x_{1}, y_{1}), G^{m}(y_{1}, x_{1})), F(F^{m}(x_{2}, y_{2}), G^{m}(y_{2}, x_{2})))\\\hspace*{5.2cm}\leq  \dfrac{k}{2}\ [d_{X}(F^{m}(x_{1},
 y_{1}), F^{m}(x_{2}, y_{2}))+ d_{Y}(G^{m}(y_{1}, x_{1}), G^{m}(y_{2}, x_{2}))]\\\hspace*{5.2cm}\leq  \dfrac{k}{2}\ \Big(\dfrac{k}{2}\Big)^{m}\ [d_{X}(x_{1}, x_{2})+d_{Y}(y_{1}, y_{2})] \\\hspace*{5.2cm}=  \Big(\dfrac{k}{2}\Big)^{m+1}\ [d_{X}(x_{1}, x_{2})+d_{Y}(y_{1}, y_{2})]$\\
Similarly,\\
$d_{Y}(G^{m+1}(y_{1}, x_{1}), G^{m+1}(y_{2}, x_{2})) \leq \Big(\dfrac{l}{2}\Big)^{m+1}\ [d_{Y}(y_{1}, y_{2})+d_{X}(x_{1}, x_{2})]$\\
Thus our claim is true for all $n\in \mathbb{N}$.\\
To prove the uniqueness we consider two cases:\\Case 1:  Assume $(x, y)$ is comparable to $(x^{*}, y^{*})$ with respect to the ordering in $X\times Y$.\\ We have,\\
$d((x, y), (x^{*}, y^{*}))= d_{X}(x, x^{*})+ d_{Y}(y, y^{*})\\\hspace*{3cm} =  d_{X}(F^{n}(x, y), F^{n}(x^{*}, y^{*}))+ d_{Y}(G^{n}(y, x), G^{n}(y^{*}, x^{*}))\\\hspace*{3cm}\leq  \Big(\dfrac{k}{2}\Big)^{n}\ [d_{X}(x, x^{*})+d_{Y}(y, y^{*})]+ \Big(\dfrac{l}{2}\Big)^{n}\ [d_{Y}(y, y^{*})+d_{X}(x, x^{*})]\\\hspace*{3cm}=  \Big\{\Big(\dfrac{k}{2}\Big)^{n}+\Big(\dfrac{l}{2}\Big)^{n}\Big\} [d_{X}(x, x^{*})+d_{Y}(y, y^{*})]\\\hspace*{3cm}\rightarrow 0 \ \ \text{as} \ \ n\rightarrow \infty $\\
This implies that $(x, y)= (x^{*}, y^{*})$.\\
Case 2: If $(x, y)$ is not comparable to $(x^{*}, y^{*})$, then by the hypothesis there exist $(u, v)\in X\times Y$ that is comparable to both $(x, y)$  and $(x^{*}, y^{*})$, which implies that $(v, u)\in Y\times X$ is comparable to both $(y, x)$  and $(y^{*}, x^{*})$ 
\\Consider\\
$d((x, y), (x^{*}, y^{*}))=  d ((F^{n}(x, y), G^{n}(y, x)), (F^{n}(x^{*}, y^{*}), G^{n}(y^{*}, x^{*})))\\\hspace*{.4cm}\leq  d((F^{n}(x, y), G^{n}(y, x)), (F^{n}(u, v), G^{n}(v, u)))\\\hspace*{.4cm} ~~+d ((F^{n}(x^{*}, y^{*}), G^{n}(y^{*}, x^{*})), (F^{n}(u, v), G^{n}(v, u)))\\\hspace*{.4cm}= d_{X}(F^{n}(x, y),
 F^{n}(u, v))+d_{Y}(G^{n}(y, x), G^{n}(v, u))+d_{X}(F^{n}(x^{*}, y^{*}), F^{n}(u, v))\\\hspace*{.4cm}~~+d_{Y}(G^{n}(y^{*}, x^{*}), G^{n}(v, u))    \\\hspace*{.4cm}\leq\Big(\dfrac{k}{2}\Big)^{n}\ [d_{X}(x, u)+d_{Y}(y, v)]+\Big(\dfrac{l}{2}\Big)^{n} \ [d_{Y}(y, v)+d_{X}(x, u)]  + \Big(\dfrac{k}{2}\Big)^{n}\ [d_{X}(x^{*}, u)+d_{Y}(y^{*}, v)]\\\hspace*{.4cm}~~+ \Big(\dfrac{l}{2}\Big)^{n}\ [d_{Y}(y^{*}, v)+d_{X}(x^{*}, u)]\\\hspace*{.4cm} = \Big\{\Big(\dfrac{k}{2}\Big)^{n}+\Big(\dfrac{l}{2}\Big)^{n}\Big\} \ [d_{X}(x, u)+d_{Y}(y, v)] +\Big\{\Big(\dfrac{k}{2}\Big)^{n}+\Big(\dfrac{l}{2}\Big)^{n}\Big\}\ [d_{X}(x^{*}, u)+d_{Y}(y^{*}, v)]\\\hspace*{.4cm}\leq \Big(\dfrac{k+l}{2}\Big)^{n}\ \big\{[d_{X}(x^{*}, u)+d_{Y}(y^{*}, v)]+[d_{X}(x, u)+d_{Y}(y, v)]\big\} \\\hspace*{.4cm}\rightarrow 0 \ \ \text{as} \ \ n\rightarrow \infty $\\
Which implies $(x, y)= (x^{*}, y^{*})$. Hence the uniqueness of FG- coupled fixed point is proved. $\square$

\begin{cly} {\upshape\lbrack 1, Theorem 2.4 \rbrack} \label{cly 3}
In addition to the hypothesis of corollary \ref{cly 1}, suppose that for all $(x, y),\ (z, t) \in X \times X$ there exists a $(u, v)\in X\times X$ that is comparable to both $(x, y)$ and $(z, t)$, then F has a unique coupled fixed point.
\end{cly}
\textbf{Proof}: Take $X= Y$, $F=G$ and $k= l$ in Theorem \ref{thm 1}, we get the result. $\square$

\begin{thm}\label{thm 5}
Let $(X, d_{X}, \leq_{P_{1}})$ and $(Y, d_{Y}, \leq_{P_{2}})$ be two partially ordered complete metric spaces. Assume that X and Y have the following properties:
\begin{enumerate}
\item[(i)] if a non decreasing sequence $\{x_{n}\} \rightarrow x$ in X, then $x_{n} \leq_{P_{1}} x$ for all n
\item[(ii)] if a non increasing sequence $\{y_{n}\} \rightarrow y$ in Y, then $y_{n} \geq_{P_{2}} y$ for all n
\end{enumerate}
Let $F: X\times Y \rightarrow X$ and $G: Y\times X\rightarrow Y$ be two functions having the mixed monotone property. Assume that there exist $k,\ l\in [0, 1)$ with 
\begin{equation*}
d_{X}(F(x, y), F(u, v))\leq \frac{k}{2}\ [d_{X}(x, u)+ d_{Y}(y, v)], \ for \ all \   x\geq_{P_{1}}u,\  \ y\leq_{P_{2}}v \tag{\ref{eqn 1}}
\end{equation*}
\begin{equation*}
d_{Y}(G(y, x), G(v, u))\leq \frac{l}{2}\ [d_{Y}(y, v)+ d_{X}(x, u)], \ for \ all \    x\leq_{P_{1}}u,\  \ y\geq_{P_{2}}v \tag{\ref{eqn 2}}
\end{equation*}
If there exist $(x_{0}, y_{0})\in X\times Y$ such that $x_{0}\leq_{P_{1}} F(x_{0}, y_{0})$ and $y_{0}\geq_{P_{2}}G(y_{0}, x_{0})$,
then there exist $(x, y)\in X\times Y$ such that $x= F(x, y)$ and $y= G(y, x)$.
\end{thm}
\textbf{Proof}:
Following the proof of Theorem \ref{thm 1} we only have to show that $(x, y)$ is an FG- coupled fixed point. Recall from the proof of Theorem \ref{thm 1} that $\{x_{n}\}$ is increasing in X and $\{y_{n}\}$ is decreasing in Y, $\lim_{n \rightarrow \infty}F^{n}(x_{0}, y_{0})= x$ and $\lim_{n \rightarrow \infty}G^{n}(y_{0}, x_{0})= y$.\\We have, \\
$d_{X}(F(x, y), x)\leq  d_{X}(F(x, y), F^{n+1}(x_{0}, y_{0}))+ d_{X}(F^{n+1}(x_{0}, y_{0}), x)\\\hspace*{2.5cm} = d_{X}(F(x, y), F(F^{n}(x_{0}, y_{0}), G^{n}(y_{0}, x_{0})))+ d_{X}(F^{n+1}(x_{0}, y_{0}), x)$\\
 By $(i)$ and $(ii)$, $x\geq_{P_{1}} F^{n}(x_{0}, y_{0})$ and $y\leq_{P_{2}} G^{n}(y_{0}, x_{0})$, therefore by (\ref{eqn 1})\\
$d_{X}(F(x, y), x)\leq \dfrac{k}{2}\ [d_{X}(x,  F^{n}(x_{0}, y_{0}))+ d_{Y}(y, G^{n}(y_{0}, x_{0})] + d_{X}(F^{n+1}(x_{0}, y_{0}), x) \\\hspace*{2.6cm} \rightarrow  0 \ \ \text{as} \ \ n\rightarrow \infty $\\Therefore we have $F(x, y)=x$.\\
Similarly we can prove that $G(y, x)= y$. This completes the proof. $\square$\\
We obtain the result of Gnana Bhaskar and Lakshmikantham \cite{2} as a corollary of our result.
 
\begin{cly} {\upshape\lbrack 1, Theorem 2.2 \rbrack}  \label{cly 2}
Let $(X, \leq)$ be a partially ordered set and suppose there is a metric d on X such that $(X, d)$ is a complete metric space. Assume that X has the following property:
\begin{enumerate}
\item[(i)] if a non decreasing sequence $\{x_{n}\} \rightarrow x$, then $x_{n} \leq x$ for all n
\item[(ii)] if a non increasing sequence $\{y_{n}\} \rightarrow y$, then $y_{n} \geq y$ for all n
\end{enumerate}
Let $F: X\times X \rightarrow X$  be a mapping having the mixed monotone property on X. Assume that there exist $k \in [0, 1)$ with 
\begin{equation*}
d(F(x, y), F(u, v))\leq \frac{k}{2} \ [d(x, u)+ d(y, v)], \ \forall \   x\geq u,\  \ y\leq v 
\end{equation*}
If there exist $x_{0}, y_{0}\in X$ such that $x_{0}\leq F(x_{0}, y_{0})$ and $y_{0}\geq F(y_{0}, x_{0})$,
then there exist $x, y \in X$ such that $x= F(x, y)$ and $y= F(y, x)$.
\end{cly}
\textbf{Proof}: Take $X= Y$, $F=G$ and $k= l$ in Theorem \ref{thm 5}, we get the result. $\square$

\begin{rk}
 \normalfont By adding to the hypothesis of Theorem \ref{thm 5} the condition: for every $(x, y), (x_{1}, y_{1})\in X\times Y$ there exists a $(u, v)\in X\times Y$ that is comparable to both $(x, y)$  and  $(x_{1}, y_{1})$, we can obtain a unique FG- coupled fixed point. 
\end{rk}

\begin{thm} \label{thm 2}
Let $(X, d_{X}, \leq_{P_{1}})$ and $(Y, d_{Y}, \leq_{P_{2}})$ be two complete partially ordered metric spaces and $F: X\times Y \rightarrow X$ and $G: Y\times X\rightarrow Y$ be two continuous functions having the mixed monotone property. Assume that there exist non negative k, l with $k+l<\ 1$ such that  
\begin{equation}\label{eqn 10}
d_{X}(F(x, y), F(u, v))\leq k\ d_{X}(x, u)+ l\ d_{Y}(y, v); \ \forall  \  x\geq_{P_{1}}u,\  \ y\leq_{P_{2}}v 
\end{equation}
\begin{equation}\label{eqn 11}
d_{Y}(G(y,  x),  G(v, u))\leq k\ d_{Y}(y,  v)+ l\ d_{X}(x, u); \ \forall \   x\leq_{P_{1}}u,\  \ y\geq_{P_{2}}v
\end{equation}
If there exist $(x_{0}, y_{0})\in X\times Y$ such that $x_{0}\leq_{P_{1}} F(x_{0}, y_{0})$ and $y_{0}\geq_{P_{2}}G(y_{0}, x_{0})$,
then there exist $(x, y)\in X\times Y$ such that $x= F(x, y)$ and $y= G(y, x)$.
\end{thm}
\textbf{Proof}: Following as in Theorem \ref{thm 1} we get an increasing sequence $\{x_{n}\}$ in X and a decreasing sequence $\{y_{n}\}$ in Y where $x_{n+1}=F(x_{n}, y_{n})=F^{n+1}(x_{0},y_{0})$ and $y_{n+1}=G(y_{n}, x_{n})=G^{n+1}(y_{0},x_{0})$.
\\ Claim: For $n \in \mathbb{N}$
\begin{equation} \label{eqn 8}
d_{X}(F^{n+1}(x_{0}, y_{0}), F^{n}(x_{0}, y_{0}))\leq (k+l)^{n}\ [d_{X}(x_{1}, x_{0})+d_{Y}(y_{1}, y_{0})]
\end{equation} 
\begin{equation} \label{eqn 9}
d_{Y}(G^{n+1}(y_{0}, x_{0}), G^{n}(y_{0}, x_{0}))\leq (k+l)^{n}\ [d_{Y}(y_{1}, y_{0})+d_{X}(x_{1}, x_{0})]
\end{equation}
By using (\ref{eqn 10}), (\ref{eqn 11}) and symmetric property of $d_{Y}$ we prove the claim.\\For $n= 1$ consider,\\
$d_{X}(F^{2}(x_{0}, y_{0}), F(x_{0}, y_{0}))= d_{X}(F(F(x_{0},y_{0}),G(y_{0},x_{0})),F(x_{0},y_{0}))\\\hspace*{4.5cm}\leq  k\ d_{X}(F(x_{0}, y_{0}), x_{0})+ l\ d_{Y}(G(y_{0},x_{0}),y_{0})\\\hspace*{4.5cm}=  k\ d_{X}(x_{1}, x_{0})+l\ d_{Y}(y_{1}, y_{0})\\\hspace*{4.5cm}\leq (k+l)\ [d_{X}(x_{1}, x_{0})+d_{Y}(y_{1}, y_{0})]$\\
Similarly,\\
$d_{Y}(G^{2}(y_{0}, x_{0}), G(y_{0}, x_{0}))\leq (k+l)\ [d_{Y}(y_{1}, y_{0})+d_{X}(x_{1}, x_{0})]$\\
Assume the result is true for $n\leq m$, then check for $n= m+1$. Consider,\\
$d_{X}(F^{m+2}(x_{0}, y_{0}), F^{m+1}(x_{0}, y_{0}))\\\hspace*{2cm}= d_{X}(F(F^{m+1}(x_{0}, y_{0}), G^{m+1}(y_{0}, x_{0})), \ F(F^{m}(x_{0}, y_{0}), G^{m}(y_{0}, x_{0})))\\\hspace*{2cm}\leq  k\ d_{X}(F^{m+1}(x_{0}, y_{0}), F^{m}(x_{0}, y_{0})) +l\ d_{Y}(G^{m+1}(y_{0}, x_{0}), G^{m}(y_{0}, x_{0}))\\\hspace*{2cm}\leq  k\ (k+l)^{m}\ [d_{X}(x_{1}, x_{0})+d_{Y}(y_{1}, y_{0})] +l\ (k+l)^{m}\ [d_{Y}(y_{1}, y_{0})+ d_{X}(x_{1}, x_{0})] \\\hspace*{2cm}\leq (k+l)^{m+1}\ [d_{X}(x_{1}, x_{0})+d_{Y}(y_{1}, y_{0})]$\\
 Similarly we can prove that\\
$d_{Y}(G^{m+2}(y_{0}, x_{0}), G^{m+1}(y_{0}, x_{0}))\leq (k+l)^{m+1}\ [d_{Y}(y_{1}, y_{0})+d_{X}(x_{1}, x_{0})]$\\
Thus the claim is true for all $n\in \mathbb{N}$.\\ Next we prove that $\{x_{n}\}$ is a Cauchy sequence in X and $\{y_{n}\}$ is a Cauchy sequence in Y using (\ref{eqn 8}) and (\ref{eqn 9}) respectively.\\ For $m\geq n$ consider,\\
$d_{X}(F^{m}(x_{0}, y_{0}), F^{n}(x_{0}, y_{0}))\\\hspace*{2cm}\leq  d_{X}(F^{m}(x_{0}, y_{0}), F^{m-1}(x_{0}, y_{0})) + d_{X}(F^{m-1}(x_{0}, y_{0}), F^{m-2}(x_{0}, y_{0}))\\\hspace*{2cm}~~+ . . . + d_{X}(F^{n+1}(x_{0}, y_{0}), F^{n}(x_{0}, y_{0}))\\\hspace*{2cm} \leq  (k+l)^{m-1}\ [d_{X}(x_{1}, x_{0})+d_{Y}(y_{1}, y_{0})] +(k+l)^{m-2}\ [d_{X}(x_{1}, x_{0})+d_{Y} (y_{1}, y_{0})]\\\hspace*{2cm}~~+ . . . + (k+l)^{n}\ [d_{X}(x_{1}, x_{0})+d_{Y}(y_{1}, y_{0})]\\\hspace*{2cm}= \{(k+l)^{m-1}+(k+l)^{m-2}+ . . . +(k+l)^{n}\} [d_{X}(x_{1}, x_{0})+d_{Y}(y_{1}, y_{0})]\\\hspace*{2cm} \leq  \dfrac{\delta^{n}}{1-\delta}\ [d_{X}(x_{1}, x_{0})+d_{Y}(y_{1}, y_{0})]; \  \text{where} \ \delta= k+l \ < \ 1\\\hspace*{2cm}\rightarrow 0 \ \text{as} \ n  \rightarrow \infty$\\
This implies that $\{F^{n}(x_{0},y_{0})\}$ is a Cauchy sequence in X. 
Similarly one can show that $\{G^{n}(y_{0}, x_{0})\}$ is a Cauchy sequence in Y. Since $(X, d_{X})$ and $(Y, d_{Y})$ are complete metric spaces we have $(x, y)\in X\times Y$ such that $\lim_{n\rightarrow \infty}F^{n}(x_{0}, y_{0})= x$ and $\lim_{n\rightarrow \infty}G^{n}(y_{0} ,x_{0})= y$. In the same lines as in Theorem \ref{thm 1} we can show that $(x, y)\in X\times Y$ is an FG- coupled fixed point. Hence the proof. $\square$
 
\begin{eg}
Let $X=(-\infty, 0]$ and $Y=[0,\infty)$ with usual order and usual metric. Define $F:X\times Y \rightarrow X$ and $G:Y \times X \rightarrow Y$ as $F(x,y)=\dfrac{4x-3y}{17}$ and $G(y,x)=\dfrac{4y-3x}{17}$, then it is easy to check that F and G satisfies the conditions $(\ref{eqn 10})$ and $(\ref{eqn 11})$ for $k=\dfrac{4}{17},\ l=\dfrac{3}{17}$. Here $(0,0)$ is the unique FG- coupled fixed point.
\end{eg}

\begin{rk}
\normalfont By adding to the hypothesis of Theorem \ref{thm 2} the condition: for every $(x, y),\ (x_{1}, y_{1})\in X\times Y$ there exists a $(u, v)\in X\times Y$ that is comparable to both $(x, y)$ and $(x_{1}, y_{1})$, we can obtain a unique FG- coupled fixed point. 
\end{rk}
In the following theorem we obtain uniqueness of FG- coupled fixed point using the above condition.

\begin{thm} \label{thm 7}
Let $(X, d_{X}, \leq_{P_{1}})$ and $(Y, d_{Y}, \leq_{P_{2}})$ be two complete partially ordered metric spaces and $F: X\times Y \rightarrow X$ and $G: Y\times X\rightarrow Y$ be two continuous functions having the mixed monotone property. Assume that for every $(x, y),\ (x_{1}, y_{1})\in X\times Y$ there exists a $(u, v)\in X\times Y$ that is comparable to both $(x, y) \ \ and \ \ (x_{1}, y_{1})$ and there exist non negative k, l with $k+l<\ 1$ such that  
\begin{equation*}
d_{X}(F(x, y), F(u, v))\leq k\ d_{X}(x, u)+ l\ d_{Y}(y, v); \ \forall  \  x\geq_{P_{1}}u,\  \ y\leq_{P_{2}}v \tag{\ref{eqn 10}}
\end{equation*}
\begin{equation*}
d_{Y}(G(y,  x),  G(v, u))\leq k\ d_{Y}(y,  v)+ l\ d_{X}(x, u); \ \forall \   x\leq_{P_{1}}u,\  \ y\geq_{P_{2}}v \tag{\ref{eqn 11}}
\end{equation*}
If there exist $(x_{0}, y_{0})\in X\times Y$ such that $x_{0}\leq_{P_{1}} F(x_{0}, y_{0})$ and $y_{0}\geq_{P_{2}}G(y_{0}, x_{0})$,
then there exist $(x, y)\in X\times Y$ such that $x= F(x, y)$ and $y= G(y, x)$.
\end{thm}
\textbf{Proof}: Following as in Theorem \ref{thm 2} we obtain existence of FG- coupled fixed point. Now we prove the uniqueness part. Suppose that $(x^{*}, y^{*})\in X\times Y$ is another FG- coupled fixed point, then we show that $d((x, y), (x^{*}, y^{*}))= 0$,\\where $x= \lim_{n\rightarrow \infty} F^{n}(x_{0}, y_{0})$ and $y= \lim_{n\rightarrow \infty}G^{n}(y_{0}, x_{0})$.\\ 
Claim: For any two points $(x_{1}, y_{1}),\ (x_{2}, y_{2})\in X\times Y$ which are comparable,
\begin{equation}\label{eqn 12}
d_{X}(F^{n}(x_{1}, y_{1}), F^{n}(x_{2}, y_{2}))\leq (k+l)^{n}\ [d_{X}(x_{1}, x_{2})+\ d_{Y}(y_{1}, y_{2})]
\end{equation}
\begin{equation}\label{eqn 13}
d_{Y}(G^{n}(y_{1}, x_{1}), G^{n}(y_{2}, x_{2}))\leq (k+l)^{n}\ [d_{Y}(y_{1}, y_{2})+\ d_{X}(x_{1}, x_{2})]
\end{equation}Without loss of generality assume that $(x_{2},y_{2})\leq(x_{1},y_{1})$.
 Using (\ref{eqn 10}) and (\ref{eqn 11}) we prove the claim.\\ For $n= 1$ we have,\\
$d_{X}(F(x_{1}, y_{1}), F(x_{2}, y_{2}))\leq k\ d_{X}(x_{1}, x_{2})+l\ d_{Y}(y_{1}, y_{2})\\ \hspace*{4.2cm}\leq (k+l)\ [d_{X}(x_{1}, x_{2})+\ d_{Y}(y_{1}, y_{2})]$\\
Now assume that the result is true for $n\leq m$ and check for $n= m+ 1$.\\Consider,\\
$d_{X}(F^{m+1}(x_{1}, y_{1}), F^{m+1}(x_{2}, y_{2}))=  d_{X}(F(F^{m}(x_{1}, y_{1}), G^{m}(y_{1}, x_{1})), F(F^{m}(x_{2}, y_{2}), G^{m}(y_{2}, x_{2})))\\\hspace*{2cm}\leq  k\ d_{X}(F^{m}(x_{1},  y_{1}), F^{m}(x_{2}, y_{2}))+l\ d_{Y}(G^{m}(y_{1}, x_{1}), G^{m}(y_{2}, x_{2}))\\\hspace*{2cm}\leq  k\ (k+l)^{m}\ [d_{X}(x_{1}, x_{2})+ d_{Y}(y_{1}, y_{2})]+l\ (k+l)^{m}\ [d_{Y}(y_{1}, y_{2})+\ d_{X}(x_{1}, x_{2})]\\\hspace*{2cm}=  (k+l)^{m+1}\ [d_{X}(x_{1}, x_{2})+\ d_{Y}(y_{1}, y_{2})]$\\Similarly we get,\\
$d_{Y}(G^{m+1}(y_{1}, x_{1}), G^{m+1}(y_{2}, x_{2}))\leq (k+l)^{m+1}\ [d_{Y}(y_{1}, y_{2})+\ d_{X}(x_{1}, x_{2})]$\\
Thus the claim is true for all $n\in \mathbb{N}$.\\
To prove the uniqueness we use the inequalities (\ref{eqn 12}) and(\ref{eqn 13}). We consider two cases:\\Case 1: Assume $(x, y)$ is comparable to $(x^{*}, y^{*})$ with respect to the ordering in $X\times Y$.\\ Now consider,\\
$d((x, y), (x^{*}, y^{*}))= d_{X}(x, x^{*})+ d_{Y}(y, y^{*})\\\hspace*{3.1cm} =  d_{X}(F^{n}(x, y),  F^{n}(x^{*}, y^{*}))+ d_{Y}(G^{n}(y, x), G^{n}(y^{*}, x^{*}))\\\hspace*{3.1cm}\leq  (k+l)^{n}\ [d_{X}(x, x^{*})+\ d_{Y}(y, y^{*})]+ (k+l)^{n}\ [d_{Y}(y, y^{*})+\ d_{X}(x, x^{*})]\\\hspace*{3.1cm}=  2\ (k+l)^{n}[d_{X}(x, x^{*})+d_{Y}(y, y^{*})]\\\hspace*{3.1cm}\rightarrow 0 \ \ \text{as} \ \ n\rightarrow \infty$\\
This implies that $(x, y)= (x^{*}, y^{*})$.\\
Case 2: If $(x, y)$ is not comparable to $(x^{*}, y^{*})$, then by the hypothesis there exist $(u, v)\in X\times Y$ that is comparable to both $(x, y)$  and $(x^{*}, y^{*})$. 
\\ Now consider\\
$d((x, y), (x^{*},y^{*}))= d((F^{n}(x, y), G^{n}(y, x)), (F^{n}(x^{*},y^{*}), G^{n}(y^{*},x^{*})))\\\hspace*{3.2cm}\leq  d((F^{n}(x, y), G^{n}(y, x)),(F^{n}(u, v), G^{n}(v, u)))\\\hspace*{3.2cm}~~+ d((F^{n}(u, v), G^{n}(v, u)),(F^{n}(x^{*},y^{*}), G^{n}(y^{*},x^{*})))\\\hspace*{3.2cm} = d_{X}(F^{n}(x, y), F^{n}(u,v))+d_{Y}(G^{n}(y, x), G^{n}(v, u)) \\\hspace*{3.2cm}~~+d_{X}(F^{n}(x^{*},y^{*}),F^{n}(u,v))+d_{Y}(G^{n}(y^{*},x^{*}),G^{n}(v, u))\\\hspace*{3.2cm}\leq (k+l)^{n}\ [d_{X}(x, u)+\ d_{Y}(y, v)]+ (k+l)^{n}\ [d_{Y}(y, v)+\ d_{X}(x, u)]\\\hspace*{3.2cm}~~+(k+l)^{n}\ [d_{X}(x^{*}, u)+\ d_{Y}(y^{*}, v)]+(k+l)^{n}\ [d_{Y}(y^{*}, v)+\ d_{X}(x^{*}, u)]\\\hspace*{3.2cm} =2\ (k+l)^{n}\ [d_{X}(x, u)+\ d_{Y}(y, v)]+ 2\ (k+l)^{n}\ [d_{X}(x^{*}, u)+\ d_{Y}(y^{*}, v)]\\\hspace*{3.2cm}\rightarrow 0 \ \ \text{as} \ \ n\rightarrow \infty$ ,\\
which implies $(x, y)= (x^{*},y^{*})$. Hence the uniqueness of FG- coupled fixed point is proved. $\square$

 The above result is valid for any two mappings F and G if the spaces satisfies a condition as shown in the following theorem.
   
\begin{thm}\label{thm 6}
Let $(X, d_{X}, \leq_{P_{1}})$ and $(Y, d_{Y}, \leq_{P_{2}})$ be two partially ordered complete metric spaces. Assume that X and Y have the following properties:
\begin{enumerate}
\item[(i)] if a non decreasing sequence $\{x_{n}\} \rightarrow x$ in X, then $x_{n} \leq_{P_{1}} x$ for all n
\item[(ii)] if a non increasing sequence $\{y_{n}\} \rightarrow y$ in Y, then $y_{n} \geq_{P_{2}} y$ for all n
\end{enumerate}
Let $F: X\times Y \rightarrow X$ and $G: Y\times X\rightarrow Y$ be two functions having the mixed monotone property. Assume that there exist non negative k, l with $k+l<1$
\begin{equation*}
d_{X} (F(x, y), F(u, v))\leq k\ d_{X}(x, u)+ l\ d_{Y}(y, v), \ \forall    x\geq_{P_{1}}u,\  \ y\leq_{P_{2}}v \tag{\ref{eqn 10}}
\end{equation*}
\begin{equation*}
d_{Y}(G(y, x), G(v, u))\leq k\ d_{Y}(y, v)+ l\ d_{X}(x, u), \ \forall    x\leq_{P_{1}}u,\  \ y\geq_{P_{2}}v \tag{\ref{eqn 11}}
\end{equation*}
If there exist $(x_{0}, y_{0})\in X\times Y$ such that $x_{0}\leq_{P_{1}} F(x_{0}, y_{0})$ and $y_{0}\geq_{P_{2}}G(y_{0}, x_{0})$,
then there exist $(x, y)\in X\times Y$ such that $x= F(x, y)$ and $y= G(y, x)$.
\end{thm}
\textbf{Proof}:
Following the proof of Theorem \ref{thm 2} we only have to show that $(x, y)$ is an FG- coupled fixed point. Recall from the proof of Theorem \ref{thm 2} that $\{x_{n}\}$ is increasing in X and $\{y_{n}\}$ is decreasing in Y, $\lim_{n \rightarrow \infty}F^{n}(x_{0}, y_{0})= x$ and $\lim_{n \rightarrow \infty}G^{n}(y_{0}, x_{0})= y$.\\We have \\
$d_{X}(F(x, y), x)\leq  d_{X}(F(x, y), F^{n+1}(x_{0}, y_{0}))+ d_{X}(F^{n+1}(x_{0}, y_{0}), x)\\\hspace*{2.7cm} = d_{X}(F(x, y), F(F^{n}(x_{0}, y_{0}), G^{n}(y_{0}, x_{0})))+ d_{X}(F^{n+1}(x_{0}, y_{0}), x)\\  \text{By (i) and (ii) we have}\ x\geq_{P_{1}} F^{n}(x_{0}, y_{0})\ and \ y\leq_{P_{2}} G^{n}(y_{0}, x_{0}).\ \text{Therefore\ using\ (\ref{eqn 10})\ we\ get} \\ d_{X}(F(x, y), x)\leq  k\ d_{X}(x, F^{n}(x_{0}, y_{0}))+ l\ d_{Y}(y, G^{n}(y_{0}, x_{0})) + d_{X}(F^{n+1}(x_{0},  y_{0}), x)\\ \hspace*{2.7cm}\rightarrow  0 \ \ \text{as} \ \ n\rightarrow \infty$\\ 
That is $F(x, y)= x$.\\Similarly we get 
$G(y, x)= y$. This completes the proof. $\square$

\begin{rk}
 \normalfont By adding to the hypothesis of Theorem \ref{thm 6} the condition: for every $(x, y), (x_{1}, y_{1})\in X\times Y$ there exists a $(u, v)\in X\times Y$ that is comparable to both $(x, y)$  and  $(x_{1}, y_{1})$, we can obtain a unique FG- coupled fixed point. 
\end{rk}

\begin{rk}
\normalfont By putting $k= l= \dfrac{k^{'}}{2}$ in theorems \ref{thm 2}, \ref{thm 7}, \ref{thm 6} we get theorems 2.1, 2.4, 2.2 of Gnana Bhaskar and Lakshmikantham \cite{2} respectively. 
\end{rk}

\begin{thm} \label{thm 3}
Let $(X, d_{X}, \leq_{P_{1}})$ and $(Y, d_{Y}, \leq_{P_{2}})$ be two partially ordered complete metric spaces and $F: X\times Y \rightarrow X$ and $G: Y\times X\rightarrow Y$ be two continuous functions having the mixed monotone property. Assume that there exist non negative k, l with $k+l<1$ such that  
\begin{equation}\label{eqn 14}
d_{X}(F(x, y), F(u, v))\leq k\ d_{X}(x, F(x, y))+ l\ d_{X}(u, F(u, v)), \ \forall    x\geq_{P_{1}}u,\  \ y\leq_{P_{2}}v 
\end{equation}
\begin{equation}\label{eqn 15}
d_{Y}(G(y, x), G(v, u))\leq k\ d_{Y}(y, G(y, x))+ l\ d_{Y}(v, G(v, u)), \ \forall    x\leq_{P_{1}}u,\  \ y\geq_{P_{2}}v
\end{equation}
If there exist $(x_{0}, y_{0})\in X\times Y$ such that $x_{0}\leq_{P_{1}} F(x_{0}, y_{0})$ and $y_{0}\geq_{P_{2}}G(y_{0}, x_{0})$,
then there exist $(x, y)\in X\times Y$ such that $x= F(x, y)$ and $y= G(y, x)$.
\end{thm}
\textbf{Proof}: 
By using the mixed monotone property of F and G and given conditions on $x_{0}$ and $y_{0}$ it is easy to show that $\{x_{n}\}$ is an increasing sequence in X and $\{y_{n}\}$ is a decreasing sequence in Y where $x_{n+1}=F(x_{n},y_{n})=F^{n+1}(x_{0},y_{0})$ and $y_{n+1}=G(y_{n},x_{n})=G^{n+1}(y_{0},x_{0})$\\ Claim: For $n \in \mathbb{N}$
\begin{equation}\label{eqn 18}
d_{X}(F^{n+1}(x_{0},y_{0}),F^{n}(x_{0},y_{0}))\leq \Big(\frac{l}{1-k}\Big)^{n} d_{X}(x_{1},x_{0})
\end{equation}  
\begin{equation}\label{eqn 19}
d_{Y}(G^{n+1}(y_{0},x_{0}),G^{n}(y_{0},x_{0}))\leq \Big(\frac{k}{1-l}\Big)^{n} d_{Y}(y_{1},y_{0})
\end{equation}Using the contraction on F and G and symmetric property on $d_{Y}$ we prove the claim.\\ For $n= 1$ consider,\\
$d_{X}(F^{2}(x_{0},y_{0}),F(x_{0},y_{0}))=  d_{X}(F(F(x_{0},y_{0}),G(y_{0},x_{0})),F(x_{0},y_{0}))\\ \hspace*{4.2cm}\leq  k\ d_{X}(F(x_{0},y_{0}),F^{2}(x_{0},y_{0}))+l\ d_{X}(x_{0},F(x_{0},y_{0}))$\\
$ \text{ie},\ (1-k)\ d_{X}(F^{2}(x_{0},y_{0}), F(x_{0},y_{0}))\leq  l\ d_{X}(x_{0},F(x_{0},y_{0}))\\\hspace*{6.4cm}= l\ d_{X}(x_{0},x_{1})\\ \text{ie},\ d_{X}(F^{2}(x_{0},y_{0}),F(x_{0},y_{0}))\leq \Big(\dfrac{l}{1-k}\Big) \ d_{X}(x_{0},x_{1})$\\
Hence for $n= 1$, the claim is true. Now assume the claim for $n\leq m$ and check for $n=m+1$. Consider,\\
$d_{X}(F^{m+2}(x_{0},y_{0}),F^{m+1}(x_{0},y_{0}))\\\hspace*{2cm}= d_{X}(F(F^{m+1}(x_{0},y_{0}),G^{m+1}(y_{0},x_{0})), F(F^{m}(x_{0},y_{0}),G^{m}(y_{0},x_{0})))\\\hspace*{2cm}\leq  k\ d_{X}(F^{m+1}(x_{0},y_{0}),F^{m+2}(x_{0},y_{0}))+l\ d_{X}(F^{m}(x_{0},y_{0}),F^{m+1}(x_{0},y_{0}))$\\
$\text{ie},\ (1-k)\ d_{X}(F^{m+2}(x_{0},y_{0}),F^{m+1}(x_{0},y_{0}))\leq l\ d_{X}(F^{m}(x_{0},y_{0}),F^{m+1}(x_{0},y_{0}))\\\hspace*{7.5cm}\leq l\ \Big(\dfrac{l}{1-k}\Big)^{m}d_{X}(x_{0},x_{1})\\  \text{ie},\ d_{X}(F^{m+2}(x_{0},y_{0}),F^{m+1}(x_{0},y_{0}))\leq \Big(\dfrac{l}{1-k}\Big)^{m+1}d_{X}(x_{0},x_{1})$\\
Similarly we get,\\
$d_{Y}(G^{m+2}(y_{0}, x_{0}), G^{m+1}(y_{0},x_{0})) \leq \Big(\dfrac{k}{1-l}\Big)^{m+1}\ d_{Y}(y_{0},y_{1})$\\
Thus our claim is true for all $n\in \mathbb{N}$. Next we prove that $\{F^{n}(x_{0},y_{0})\}$ and $\{G^{n}(y_{0},x_{0})\}$ are Cauchy sequences in X and Y using (\ref{eqn 18}) and (\ref{eqn 19}) respectively.\\ For $m>n$ consider,\\
$d_{X}(F^{m}(x_{0},y_{0}),F^{n}(x_{0},y_{0}))\leq d_{X}(F^{m}(x_{0},y_{0}),F^{m-1}(x_{0},y_{0}))\\\hspace*{1cm}~~+d_{X}(F^{m-1}(x_{0},y_{0}),F^{m-2}(x_{0},y_{0})) + . . .+d_{X}(F^{n+1}(x_{0},y_{0}),F^{n}(x_{0},y_{0}))\\\hspace*{1cm}\leq \Big(\dfrac{l}{1-k}\Big)^{m-1}\ d_{X}(x_{0},x_{1})+\Big(\dfrac{l}{1-k}\Big)^{m-2}\ d_{X}(x_{0},x_{1}) + . . . +\Big(\dfrac{l}{1-k}\Big)^{n}\ d_{X}(x_{0},x_{1})\\\hspace*{1cm} = \Big\{\Big(\dfrac{l}{1-k}\Big)^{m-1}+\Big(\dfrac{l}{1-k}\Big)^{m-2}+ . . . +\Big(\dfrac{l}{1-k}\Big)^{n}\Big\}\ d_{X}(x_{0},x_{1})\\\hspace*{1cm}\leq  \Big(\dfrac{\delta^{n}}{1-\delta}\Big)\ d_{X}(x_{0},x_{1}) ; \ \text{where} \ \delta  = \ \dfrac{l}{1-k} \ <  \ 1\\\hspace*{1cm}   \rightarrow 0 \ \text{ as}  \ n\rightarrow \ \infty$\\
This implies that $\{F^{n}(x_{0},y_{0})\}$ is a Cauchy sequence in X. In a similar manner we can prove that $\{G^{n}(y_{0},x_{0})\}$ is a Cauchy sequence in Y.\\
Since $(X,d_{X})$ and $(Y,d_{Y})$ are complete metric spaces we have $(x,y)\in X\times Y$ such that 
$\lim_{n \rightarrow \infty}F^{n}(x_{0},y_{0})= x$ and $\lim_{n \rightarrow \infty}G^{n}(y_{0},x_{0})= y$. Proceeding as in Theorem \ref{thm 1} we get $(x,y)$ is an FG- coupled fixed point. This completes the proof. $\square$

\begin{eg}
Let $X=[1,2]$ and $Y=[-2,-1]$ with usual metric and usual order. Define $F:X\times Y \rightarrow X$ and $G:Y \times X \rightarrow Y$ by $F(x,y)=\dfrac{x}{4}+1$ and $G(y,x)=\dfrac{y}{4}-1$, then we can see that the conditions $(\ref{eqn 14})$ and $(\ref{eqn 15})$ for F and G  are satisfied with $k=\dfrac{1}{3},\ l=\dfrac{1}{2}$. Here       $(\dfrac{4}{3},\ \dfrac{-4}{3} )$ is the FG- coupled fixed point.
\end{eg}

\begin{thm}
Let $(X, d_{X}, \leq_{P_{1}})$ and $(Y, d_{Y}, \leq_{P_{2}})$ be two partially ordered complete metric spaces. Assume that X and Y have the following properties:
\begin{enumerate}
\item[(i)] if a non decreasing sequence $\{x_{n}\} \rightarrow x$ in X, then $x_{n} \leq_{P_{1}} x$ for all n
\item[(ii)] if a non increasing sequence $\{y_{n}\} \rightarrow y$ in Y, then $y_{n} \geq_{P_{2}} y$ for all n
\end{enumerate}
Let $F: X\times Y \rightarrow X$ and $G: Y\times X\rightarrow Y$ be two functions having the mixed monotone property. Assume that there exist non negative k, l with $k+l<1$ such that
\begin{equation*}
d_{X}(F(x, y), F(u, v))\leq k\ d_{X}(x, F(x, y))+ l\ d_{X}(u, F(u, v)), \ \forall    x\geq_{P_{1}}u,\  \ y\leq_{P_{2}}v \tag{\ref{eqn 14}}
\end{equation*}
\begin{equation*}
d_{Y}(G(y, x), G(v, u))\leq k\ d_{Y}(y, G(y, x))+ l\ d_{Y}(v, G(v, u)), \ \forall    x\leq_{P_{1}}u,\  \ y\geq_{P_{2}}v \tag{\ref{eqn 15}}
\end{equation*}
If there exist $(x_{0}, y_{0})\in X\times Y$ such that $x_{0}\leq_{P_{1}} F(x_{0}, y_{0})$ and $y_{0}\geq_{P_{2}}G(y_{0}, x_{0})$,
then there exist $(x, y)\in X\times Y$ such that $x= F(x, y)$ and $y= G(y, x)$.
\end{thm}
\textbf{Proof}:
Following as in the proof of Theorem \ref{thm 3} it remains to show that $(x, y)$ is an FG- coupled fixed point. Recall from the proof of Theorem \ref{thm 3} that $\{x_{n}\}$ is increasing in X and $\{y_{n}\}$ is decreasing in Y, $\lim_{n \rightarrow \infty}F^{n}(x_{0},y_{0})= x$ and $\lim_{n \rightarrow \infty}G^{n}(y_{0},x_{0})= y$.\\We have\\ 
$d_{X}(F(x, y),x)\leq  d_{X}(F(x, y), F^{n+1}(x_{0}, y_{0}))+ d_{X}(F^{n+1}(x_{0}, y_{0}),x)\\\hspace*{3cm} = d_{X}(F(x, y), F(F^{n}(x_{0},y_{0}),G^{n}(y_{0},x_{0})))+ d_{X}(F^{n+1}(x_{0}, y_{0}),x)\\  \text{By (i) and (ii) we have}\ x\geq_{P_{1}} F^{n}(x_{0}, y_{0})\ and \ y\leq_{P_{2}} G^{n}(y_{0}, x_{0}).\ \text{Therefore}  \ \text{using (\ref{eqn 14})\ we\ get}\\ d_{X}(F(x, y),x)\leq  k\ d_{X}(x, F(x, y))+ l\ d_{X}(F^{n}(x_{0},y_{0}), F^{n+1}(x_{0},y_{0}))\\\hspace*{3cm}+ d_{X}(F^{n+1}(x_{0}, y_{0}),x)\\ As \ n \rightarrow \infty, \ d_{X}(F(x, y),x)\leq k\ d_{X}(x, F(x, y))$. \\ 
This is possible if $d_{X}(F(x, y),x)=0$. Hence we have $F(x, y)= x$.\\Similarly using (\ref{eqn 15}) we  get 
$d_{Y}(G(y, x),y)\leq   l\ d_{Y}(y, G(y,x))$.\\
This is possible if $d_{Y}(y, G(y,x))= 0$. Hence $G(y, x)= y$. This completes the proof. $\square$ 

\begin{thm} \label{thm 4}
Let $(X, d_{X}, \leq_{P_{1}})$ and $(Y, d_{Y}, \leq_{P_{2}})$ be two partially ordered complete metric spaces and $F: X\times Y \rightarrow X$ and $G: Y\times X\rightarrow Y$ be two continuous functions having the mixed monotone property. Assume that there exist $k, l\in [0,\frac{1}{2})$  such that  
\begin{equation}\label{eqn 16}
d_{X}(F(x, y), F(u, v))\leq k\ d_{X}(x, F(u, v))+ l\ d_{X}(u, F(x, y)), \ \forall    x\geq_{P_{1}}u,\  \ y\leq_{P_{2}}v 
\end{equation}
\begin{equation}\label{eqn 17}
d_{Y}(G(y, x), G(v, u))\leq k\ d_{Y}(y, G(v, u))+ l\ d_{Y}(v, G(y, x)), \ \forall    x\leq_{P_{1}}u,\  \ y\geq_{P_{2}}v
\end{equation}
If there exist $(x_{0}, y_{0})\in X\times Y$ such that $x_{0}\leq_{P_{1}} F(x_{0}, y_{0})$ and $y_{0}\geq_{P_{2}}G(y_{0}, x_{0})$,
then there exist $(x, y)\in X\times Y$ such that $x= F(x, y)$ and $y= G(y, x)$.
\end{thm}
\textbf{Proof}: 
As in Theorem \ref{thm 1} we can construct an increasing sequence $\{x_{n}\}$ in X and a decreasing sequence $\{y_{n}\}$ in Y where $x_{n+1}=F(x_{n},y_{n})=F^{n+1}(x_{0},y_{0})$ and $y_{n+1}=G(y_{n},x_{n})=G^{n+1}(y_{0},x_{0})$.\\ Claim: For $n \in \mathbb{N}$
\begin{equation}\label{eqn 20}
d_{X}(F^{n+1}(x_{0},y_{0}),F^{n}(x_{0},y_{0}))\leq \Big(\frac{l}{1-l}\Big)^{n}\ d_{X}(x_{1},x_{0})
\end{equation} 
\begin{equation}\label{eqn 21}
d_{Y}(G^{n+1}(y_{0},x_{0}),G^{n}(y_{0},x_{0}))\leq \Big(\frac{k}{1-k}\Big)^{n}\ d_{Y}(y_{1},y_{0})
\end{equation} Using (\ref{eqn 16}), (\ref{eqn 17}) and symmetric property of $d_{Y}$ we prove the claim.\\ For $n= 1$, consider,\\
$d_{X}(F^{2}(x_{0},y_{0}),F(x_{0},y_{0}))=  d_{X}(F(F(x_{0},y_{0}),G(y_{0},x_{0})),F(x_{0},y_{0}))\\\hspace*{4.6cm} \leq  k\ d_{X}(F(x_{0},y_{0}),F(x_{0},y_{0}))+l\ d_{X}(x_{0},F^{2}(x_{0},y_{0}))\\ \hspace*{4.6cm}\leq  l\ [d_{X}(x_{0},F(x_{0},y_{0}))+d_{X}(F(x_{0},y_{0}),F^{2}(x_{0},y_{0}))]$ \\
$\text{ie},\  (1-l)\ d_{X}(F^{2}(x_{0},y_{0}), F(x_{0},y_{0}))\leq l\ d_{X}(x_{0},F(x_{0},y_{0}))\\\hspace*{6.3cm}= l\ d_{X}(x_{0},x_{1})\\ \text{ie},\  d_{X}(F^{2}(x_{0},y_{0}),F(x_{0},y_{0}))\leq \Big(\dfrac{l}{1-l}\Big) \ d_{X}(x_{0},x_{1})$\\
ie, for $n= 1$, the claim is true.\\ Now assume the claim for $n\leq m$ and check for $n=m+1$. Consider,\\
$d_{X}(F^{m+2}(x_{0},y_{0}),F^{m+1}(x_{0},y_{0}))$\\\hspace*{2cm}$= d_{X}(F(F^{m+1}(x_{0},y_{0}),G^{m+1}(y_{0},x_{0})), F(F^{m}(x_{0},y_{0}),G^{m}(y_{0},x_{0})))$\\\hspace*{2cm}$\leq  k\ d_{X}(F^{m+1}(x_{0},y_{0}),F^{m+1}(x_{0},y_{0}))+l\ d_{X}(F^{m}(x_{0},y_{0}),F^{m+2}(x_{0},y_{0}))$\\ \hspace*{2cm}$\leq l\ [d_{X}(F^{m}(x_{0},y_{0}),F^{m+1}(x_{0},y_{0}))+  d_{X}(F^{m+1}(x_{0},y_{0}),F^{m+2}(x_{0},y_{0}))]$\\
ie, $ (1-l)\ d_{X}(F^{m+2}(x_{0},y_{0}),F^{m+1}(x_{0},y_{0}))\leq l\ d_{X}(F^{m}(x_{0},y_{0}),F^{m+1}(x_{0},y_{0}))$\\$\hspace*{7.2cm}\leq l\ \Big(\dfrac{l}{1-l}\Big)^{m}\ d_{X}(x_{0},x_{1})$\\  ie, $ d_{X}(F^{m+2}(x_{0},y_{0}),F^{m+1}(x_{0},y_{0}))\leq \Big(\dfrac{l}{1-l}\Big)^{m+1}\ d_{X}(x_{0},x_{1})$\\
Similarly we get
$d_{Y}(G^{m+2}(y_{0}, x_{0}), G^{m+1}(y_{0},x_{0}))\leq \Big(\dfrac{k}{1-k}\Big)^{m+1}\ d_{Y}(y_{0},y_{1})$\\ Thus the claim is true for all $n\in \mathbb{N}$. Now using (\ref{eqn 20}) and (\ref{eqn 21}) we  prove that $\{F^{n}(x_{0},y_{0})\}$ and $\{G^{n}(y_{0},x_{0})\}$ are Cauchy sequences in X and Y respectively.\\ For $m>n$, consider,\\
$d_{X}(F^{m}(x_{0},y_{0}),F^{n}(x_{0},y_{0}))$ \\\hspace*{2cm} $\leq d_{X}(F^{m}(x_{0},y_{0}),F^{m-1}(x_{0},y_{0}))+d_{X}(F^{m-1}(x_{0},y_{0}),F^{m-2}(x_{0},y_{0}))$\\\hspace*{2cm}$~~+ . . .+d_{X}(F^{n+1}(x_{0},y_{0}),F^{n}(x_{0},y_{0}))$\\\hspace*{2cm}$\leq \Big(\dfrac{l}{1-l}\Big)^{m-1}\ d_{X}(x_{0},x_{1})+\Big(\dfrac{l}{1-l}\Big)^{m-2}\ d_{X}(x_{0},x_{1})$\\ \hspace*{2cm}$~~+ . . . +\Big(\dfrac{l}{1-l}\Big)^{n}\ d_{X}(x_{0},x_{1})$\\\hspace*{2cm}$ = \Big\{\Big(\dfrac{l}{1-l}\Big)^{m-1}+\Big(\dfrac{l}{1-l}\Big)^{m-2}+ . . . +\Big(\dfrac{l}{1-l}\Big)^{n}\Big\}\ d_{X}(x_{0},x_{1})$\\ \hspace*{2cm}$\leq  \Big(\dfrac{\delta^{n}}{1-\delta}\Big)\ d_{X}(x_{0},x_{1}) ; \ \text{where} \ \delta  = \ \dfrac{l}{1-l} \ <  \ 1$\\  \hspace*{2cm} $\rightarrow 0 \ \text{as}  \ n\rightarrow \ \infty$\\
This implies that $\{F^{n}(x_{0},y_{0})\}$ is a Cauchy sequence in X.\\Similarly we prove that $\{G^{n}(y_{0},x_{0})\}$ is a Cauchy sequence in Y.
Since $(X,d_{X})$ and $(Y,d_{Y})$ are complete metric spaces we have $(x,y)\in X\times Y$ such that 
$\lim_{n \rightarrow \infty}F^{n}(x_{0},y_{0})= x$ and $\lim_{n \rightarrow \infty}G^{n}(y_{0},x_{0})= y$. \\By using the continuity of F and G we can prove that $(x,y)$ is an FG- coupled fixed point as in the Theorem \ref{thm 1}. Hence the result. $\square$

\begin{eg}
Let $X=[0,1]$ with usual metric, $x,u\in [0,1]$ with $x\leq_{P_{1}}u \ \Leftrightarrow \  x=u \ $and $Y=[-1,0]$ with usual metric, $y,v\in [-1,0]$ with $y\leq_{P_{2}}v \ \Leftrightarrow \ \text{either}\ y=v \ \text{or}\ (y,v)=(-1,0)$. Define $F:X\times Y \rightarrow X$ and $G:Y \times X \rightarrow Y$ as $F(x,y)=\dfrac{x}{3}$ and $G(y,x)=\dfrac{-x}{3}$, then we can see that the conditions $(\ref{eqn 16})$ and $(\ref{eqn 17})$ for F and G  are satisfied for any $k,\ l\in [0,\frac{1}{2})$. Here $(0,0)$ is the unique FG- coupled fixed point.
\end{eg}
We can replace the continuity of F and G by other conditions to obtain FG- coupled fixed point result as follows:

\begin{thm}
Let $(X, d_{X}, \leq_{P_{1}})$ and $(Y, d_{Y}, \leq_{P_{2}})$ be two partially ordered complete metric spaces. Assume that X and Y have the following properties:
\begin{enumerate}
\item[(i)] if a non decreasing sequence $\{x_{n}\} \rightarrow x$ in X, then $x_{n} \leq_{P_{1}} x$ for all n
\item[(ii)] if a non increasing sequence $\{y_{n}\} \rightarrow y$ in Y, then $y_{n} \geq_{P_{2}} y$ for all n
\end{enumerate}
Let $F: X\times Y \rightarrow X$ and $G: Y\times X\rightarrow Y$ be two functions having the mixed monotone property. Assume that there exist $k, l\in [0, \frac{1}{2})$ with 
\begin{equation*}
d_{X}(F(x, y), F(u, v))\leq k\ d_{X}(x, F(u, v))+ l\ d_{X}(u, F(x, y)), \ \forall    x\geq_{P_{1}}u,\  \ y\leq_{P_{2}}v \tag{\ref{eqn 16}}
\end{equation*}
\begin{equation*}
d_{Y}(G(y, x), G(v, u))\leq k\ d_{Y}(y, G(v, u))+ l\ d_{Y}(v, G(y, x)), \ \forall    x\leq_{P_{1}}u,\  \ y\geq_{P_{2}}v \tag{\ref{eqn 17}}
\end{equation*}
If there exist $(x_{0}, y_{0})\in X\times Y$ such that $x_{0}\leq_{P_{1}} F(x_{0}, y_{0})$ and $y_{0}\geq_{P_{2}}G(y_{0}, x_{0})$,
then there exist $(x, y)\in X\times Y$ such that $x= F(x, y)$ and $y= G(y, x)$.
\end{thm}
\textbf{Proof}:
Following as in the proof of Theorem \ref{thm 4} we only have to show that $(x, y)$ is an FG- coupled fixed point. Recall from the proof of Theorem \ref{thm 4} that $\{x_{n}\}$ is increasing in X and $\{y_{n}\}$ is decreasing in Y, $\lim_{n \rightarrow \infty}F^{n}(x_{0},y_{0})= x$ and $\lim_{n \rightarrow \infty}G^{n}(y_{0},x_{0})= y$.\\Now consider \\
$d_{X}(F(x, y),x)\leq  d_{X}(F(x, y), F^{n+1}(x_{0}, y_{0}))+ d_{X}(F^{n+1}(x_{0}, y_{0}),x)\\ \hspace*{2.5cm}= d_{X}(F(x, y), F(F^{n}(x_{0},y_{0}),G^{n}(y_{0},x_{0})))+ d_{X}(F^{n+1}(x_{0}, y_{0}),x)\\   \text{By (i) and (ii) we have}\ x\geq_{P_{1}} F^{n}(x_{0}, y_{0})\ and \ y\leq_{P_{2}} G^{n}(y_{0}, x_{0}).\ \text{Therefore\ using\ (\ref{eqn 16}) we\ get} \\ d_{X}(F(x, y),x)\leq  k\ d_{X}(x, F^{n+1}(x, y))+ l\ d_{X}(F^{n}(x_{0},y_{0}), F(x,y))+ d_{X}(F^{n+1}(x_{0}, y_{0}),x)\\ As \ n \rightarrow \infty, \ d_{X}(F(x, y),x)\leq l\ d_{X}(x, F(x, y))$ \\ 
This is possible if $d_{X}(F(x, y),x)=0$. Hence we have $F(x, y)= x$.\\
Similarly using (\ref{eqn 17}) we prove that $d_{Y}(y, G(y,x))= 0$. Thus $G(y, x)= y$. This completes the proof. $\square$

\begin{rk}
 \normalfont In all the theorems in this paper if we put $X= Y$ and $F= G$ we get several coupled fixed point theorems in partially ordered complete metric spaces.
\end{rk}

\section*{Acknowledgement}
The first author would like to thank Kerala State Council for Science, Technology and Environment for the financial support.

\hspace*{-.7cm}\textsc{ Prajisha Eacha} (Corresponding author)\\
\textsc{Department of Mathematics\\
Central University of Kerala, India}\\
\textit{E-mail address: prajisha1991@gmail.com}\\

\hspace*{-.7cm}\textsc{ Shaini Pulickakunnel\\
Department of Mathematics\\
Central University of Kerala, India}\\
\textit{E-mail address: shainipv@gmail.com}

\end{document}